\documentclass[10pt]{article}

\usepackage{bm}
\usepackage[top=0.9in,bottom=1.1in,left=1.05in,right=1.05in]{geometry}

\usepackage{amsmath}
\usepackage{graphicx}
\usepackage[colorinlistoftodos]{todonotes}
\usepackage[colorlinks=true, allcolors=blue]{hyperref}
\usepackage{amsmath}
\usepackage{upgreek}
\usepackage{amssymb}
\usepackage{amsfonts}
\usepackage{amsthm}
\usepackage{mathrsfs}
\usepackage{fontenc}
\usepackage{empheq}
\usepackage{enumerate}
\usepackage{comment}
\usepackage{appendix}
\usepackage{mathtools}
\usepackage{booktabs}
\mathtoolsset{showonlyrefs}

\hypersetup{urlcolor=black}

\theoremstyle{plain}
\newtheorem{thm}{Theorem}

\newtheorem{prop}[thm]{Proposition}

\newcommand{\RR}{\mathbb{R}}
\newcommand{\EE}{\mathbb{E}}

\newcommand{\PP}{\mathbb{P}}

\newcommand{\MCF}{\mathcal{F}}
\newcommand{\ud}{\,\mathrm{d}}
\newcommand{\mc}[1]{\mathcal{#1}}
\newcommand{\half}{\frac{1}{2}}

\newcommand{\abs}[1]{\left|#1\right|}

\theoremstyle{definition}

\theoremstyle{remark}
\newtheorem{rem}{Remark}

\newcommand{\transpose}{^{\operatorname{T}}}
\newcommand{\inv}{^{\raisebox{.2ex}{$\scriptscriptstyle-1$}}}

\title{Recurrent Neural Networks for Stochastic Control Problems with Delay}
\author{Jiequn Han\thanks{Department of Mathematics, Princeton University, Princeton, NJ 08544-1000, USA, \em{jiequnh@princeton.edu}.} \and Ruimeng Hu\thanks{Department of Mathematics, and Department of Statistics and Applied Probability, University of California, Santa Barbara, CA 93106-3080, USA, {\em rhu@ucsb.edu}. RH was partially supported by the NSF grant DMS-1953035, and the Faculty Career Development Award and the Research Assistant Program Award at UCSB.}}
\date{\today}

\begin{document}

\maketitle
\begin{abstract}

Stochastic control problems with delay are challenging due to the path-dependent feature of the system and thus its intrinsic high dimensions. In this paper, we propose and systematically study deep neural networks-based algorithms to solve stochastic control problems with delay features. Specifically, we employ neural networks for sequence modeling (\emph{e.g.}, recurrent neural networks such as long short-term memory) to parameterize the policy and optimize the objective function. The proposed algorithms are tested on three benchmark examples: a linear-quadratic problem, optimal consumption with fixed finite delay, and portfolio optimization with complete memory. Particularly, we notice that the architecture of recurrent neural networks naturally captures the path-dependent feature with much flexibility and yields better performance with more efficient and stable training of the network compared to feedforward networks. The superiority is even evident in the case of portfolio optimization with complete memory, which features infinite delay.

\end{abstract}

\textbf{Keywords:}
Deep learning, stochastic control with delay, recurrent neural networks, stochastic differential delay equations.

\section{Introduction}

Stochastic control problems study the agent's rational behavior with the existence of uncertainty in observations or in the noise that drives the evolution of the system. Inclusion of delay in stochastic control problems is important for realistic applications, {\it e.g.}, in economics for time-to-build problems \cite{kydland1982time,asea1999time}, in marketing for modeling the ``carryover'' or ``distributed lag''   advertising effect \cite{gozzi200513,gozzi2009controlled}, and in finance for portfolio selection under the market with memory and delayed responses  \cite{oksendal2000maximum,federico2011stochastic,elsanosi2001optimal,li2018portfolio}. See also \cite[Chapter 1]{kolmanovskiui1996control} for modeling systems with aftereffect in mechanics and engineering, biology, and medicine. To model the delay feature, the dynamics of the controlled system will depend not only on the current state but also on the history of $\delta$ time units prior to the current time, where $\delta$ is a fixed number and can be infinity. This makes the problem path-dependent and, thus, infinite-dimensional.

The challenge brought by the path-dependence feature in stochastic control problems with delay has attracted rich theoretical studies in the literature. For example, \cite{gozzi200513,gozzi2009controlled,gozzi2017stochastic} proposed to reformulate them as infinite-dimensional Markovian problems and analyze the associated Hamilton Jacobi Bellman (HJB) equations. The dynamic programming principle was proved in \cite{larssen2002dynamic,chen2012dynamic}. In \cite{peng2009anticipated}, the authors studied the problem using anticipated backward stochastic differential equations (ABSDEs), and \cite{chen2010maximum,guatteri2020stochastic,oksendal2011optimal} proved a stochastic maximum principle using the ABSDEs. \cite{bandini2018backward} employed the so-called randomization method, to list a few. Meanwhile, except for the special cases where problems can be reduced to finite-dimensional ones \cite{elsanosi2000some,larssen2001hjb,oksendal2000maximum,elsanosi2001optimal,bauer2005stochastic}, the stochastic control problems with delay remain practically intractable, and one needs to resort to numerical methods for possible solutions; see \cite{kushner2006numerical,fischer2007discretisation,fischer2008time,kushner2008numerical} for probabilistic approaches which analyzed the corresponding discretized control problems, and  \cite{chang2008finite} for a analytical approach which focused on finite difference methods of the infinite-dimensional HJB equation. In both approaches, one has to temporally and spatially discretize the stochastic control problems with delay, yielding a finite-dimensional setup with dimensionality proportional to both the number of discretized timestamps and spatial grids. Therefore, all aforementioned existing algorithms can only work in the low-dimensional setting but encounter demanding challenges or become unfeasible  when faced with high-dimensional cases.

This paper aims to address the aforementioned numerical challenges by deep learning-based algorithms, with the model described by a stochastic differential delay equation (SDDE). 
Observing deep neural networks' remarkable performance in representing high-dimensional functions in numerical computations in many fields \cite{Be:09,CaTr:17,EHaJe:17,HaJeE:18,han2019uniformly,HaZhE:19,han2020solving,pfau2020ab,hermann2020deep}, we naturally leverage them in the context of stochastic control problems with delay. 
Specifically, motivated by \cite{HaE:16}, we shall approximate the controls using neural networks of various architectures at each time, stack these subnetworks together to form a deep network, and train them simultaneously. The optimal parameters are obtained by minimizing the loss function, which is the proxy of the cost functional in the control problem. Although using deep neural networks to direct parametrize the strategy in optimal control problems is not new, {\it e.g.}, in \cite{HaE:16,Hu2:19,carmona2019convergence,fouque2020deep}, our work has the following merits: Firstly, we develop deep learning algorithms with the focus on the feature of delay, which is by nature infinite-dimensional; To the best of the authors' knowledge, this is the first work in the literature that systematically leverages neural networks to solve stochastic control problems with delay beyond the linear-quadratic case.
Secondly, we systematically study the strengths and weaknesses of different neural network architecture by testing them on three typical examples with benchmark solutions. The carefully selected benchmark problems with open-sourced code facilitate the further study of numerical algorithms for stochastic control problems with delay features. Our main findings include: 
\begin{itemize}
    \item[(1)]The algorithm based on recurrent neural networks can naturally capture the  path-dependent feature, \emph{i.e.}, not requiring a priori knowledge of the lag time $\delta$, and yield better performance with more efficient and stable training compared to feedforward neural networks; 
    \item[(2)] The former one is capable of dealing with more complex problems efficiently and accurately, {\it e.g.}, problems with infinite delay $\delta = \infty$, and problems with state constraints;
    \item[(3)]  Both algorithms perform better when training with state processes (corresponds to closed-loop controls) than with background noise (corresponds to open-loop controls), especially for problems with state constraints.
\end{itemize}

The rest of the paper is organized as follows. In Section~\ref{sec_formulation}, we introduce the mathematical formulation of the stochastic control problems with delay in continuous time. We describe the deep learning-based algorithms in Section~\ref{sec_algorithm} and three benchmark examples in Section~\ref{sec_benchmark}, followed by a systemic numerical study in Section~\ref{sec_numerics}. We make conclusive remarks and describe future works in Section~\ref{sec_conclusion}.

\section{The stochastic control problem with delay}\label{sec_formulation}
On a complete probability space $(\Omega, \mc{F}, \PP)$, we consider a stochastic control problem in which the state process $X \in \RR^n$ is characterized by a stochastic  differential delay equation (SDDE):
\begin{equation}\label{def:Xt}
\begin{dcases}
    \ud X(t) = b(t, X_t, \pi(t)) \ud t + \sigma(t, X_t, \pi(t)) \ud W(t),  & t \in [0, T],\\
    X(t) = \varphi(t), & t \in [-\delta, 0].
\end{dcases}
\end{equation}
Here $\delta \geq 0$ is the fixed delay, $\pi$ is the control process taking values in $\mc{A} \in \RR^m$  and to be chosen in some admissible set $\mathbb{A}$ and $W(t)$ is an $\ell$-dimensional standard Brownation motion. 
Throughout the paper, we denote by $P_t$ the trajectory of a process $P$ from time $t-\delta$ to $t$, and $P(t)$ the process value at time $t$, \emph{i.e.}, $P_t(s) = P(t+ s)$, for $-\delta \leq s \leq 0$. 

\begin{rem}
	The delay parameter $\delta$ here is assumed to be deterministic and known. However, in a more realistic setting $\delta$ might be \emph{unknown} a priori. We shall see that the recurrent neural networks (discussed in Section~\ref{sec:lstm}) can naturally take care of the unknown $\delta$ implicitly due to its architecture design. 
\end{rem}

Let $C := C([-\delta, 0], \RR^n)$ be the Banach space of all continuous functions with the supremum norm:
\begin{equation}
   \|y\|_C = \sup_{-\delta \leq s \leq 0}|y(s)|, \;\forall y \in C.
\end{equation}
The drift $b$ and volatility $\sigma$ coefficients are deterministic functionals:
\begin{equation}
    (b,\sigma) : [0,T] \times C \times \mc{A} \to (\RR^n, \RR^{n \times \ell}).
\end{equation}
Denote by $L^2(\Omega, C)$ the space of all $\MCF$-measurable stochastic processes, \emph{i.e.},
\begin{equation}
    \Omega \ni \omega \to X(\omega) \in C \text{ is in } L^2(\Omega, C), \text{ iff. } \int_\Omega \|X(\omega)\|_C^2 \ud \PP(\omega) < \infty,
\end{equation}
then $L^2(\Omega, C)$ is complete with the semi-norm $\|X\|_{L^2(\Omega, C)} := [\int_\Omega \|X(\omega)\|_C^2 \ud \PP(\omega)]^{1/2}$. We assume that the initial path $\varphi \in L^2(\Omega, C)$ and is independent of the Brownian motion $W(t)$,  and the existence of solution $X$ to the SDDE \eqref{def:Xt} is considered in $L^2(\Omega, C([-\delta, T], \RR^n))$. Let $(\mc{F}_t)_{t\geq 0}$ be the filtration supporting $W(t)$ and $\varphi$, and let $C([0,T], L^2(\Omega, C))$ be the space of all $L^2$-continuous $C$-valued $\MCF_t$-adapted process $P: [0,T] \ni t \to P_t \in L^2(\Omega, C)$ with the semi-norm:
\begin{equation}
    \|P\|_{C([0,T], L^2(\Omega, C))} := \sup_{0 \leq t \leq T} \| P_t \|_{L^2(\Omega, C)}.
\end{equation}
The trajectory $X_t$ of SDDE \eqref{def:Xt} is considered in $C([0,T], L^2(\Omega, C))$.

The agent aims to minimize her expected cost:
\begin{equation}\label{def:cost}
\EE_{\varphi}\left[\int_0^T f(t, X_t, \pi(t)) \ud t  + g(X_T)\right],
\end{equation}
for a given distribution of the initial condition $\varphi$ and over all admissible strategies $\pi$ in $\mathbb{A}$:
\begin{equation}\label{def_A}
\begin{aligned}
    \mathbb{A} := &\left\{\{\mc{F}_t\}\text{-progressively measurable process } \pi: [0,T] \times \Omega \to \mc{A} \subset \RR^m:
  \int_0^T \EE[\pi(t)^2]\ud t <\infty \right\}.
\end{aligned}
\end{equation}
where the running cost $f$ and terminal cost $g$ are deterministic functionals, $f: [0,T] \times C \times \mc{A} \to \RR$,  $g: C \to \RR$.

Usually, one requires uniform Lipschitz conditions in the second variable of $b$ and $\sigma$ to ensure the existence and uniqueness of strong solutions to SDDE \eqref{def:Xt}, that is,
\begin{equation}
    \|(b, \sigma)(t, y_1, \pi) - (b, \sigma)(t, y_2, \pi) \|_{L^2} \leq L \|y_1 - y_2\|_{L^2(\Omega, C)}, \; \forall t \in [0, T] \text{ and } y_1, y_2 \in L^2(\Omega, C). 
\end{equation}
See detailed analysis in Mohammed's monographs \cite{mohammed1984stochastic,mohammed1998stochastic}. Assumptions on $f$ and $g$ would ensure the expected cost \eqref{def:cost} is finite. 

In this paper, instead of discussing necessary conditions for the admissibility, we aim at providing a systematic numerical study of deep learning algorithms for finding optimal strategies to stochastic control problems with delay. That is, we focus on the deep neural networks' (DNNs) architecture design in order to handle the high-dimensionality arising from the delay and comparing their performance based on some tractable examples. 
In Section~\ref{sec_benchmark}, we present three examples with tractability to benchmark our numerical schemes and support our findings. Note that the third example has an infinite history dependence (\emph{i.e.}, $\delta = \infty$), and the numerical results further illustrate the advantage of using recurrent neural networks in this more general frameworks.

\section{Deep learning algorithm}\label{sec_algorithm}
Our numerical algorithm builds on the temporal discretization of \eqref{def:Xt}--\eqref{def:cost} and approximating  $\pi(t) \in \RR^m$ using neural networks. More precisely, let $N_T \in \mathbb{N}$ and $0 = t_0 < t_1 < \ldots < t_{N_T} = T$ be a partition of size $N_T$ on $[0,T]$.  Without loss of generality, we assume they are equidistributed and the fixed delay $\delta <\infty$ covers $N_\delta$ subintervals:
\begin{equation}
  h \equiv t_{k+1} - t_k, \; \forall k = 0,  \ldots,  N_T-1, \text{ and } \delta = N_\delta h.
\end{equation}
Consequently we can extend the partition to $[-\delta, 0]$:
\begin{equation}
 \; -\delta = t_{-N_\delta} \leq t_{-N_\delta+1} \leq \ldots t_0 = 0, \text{ with } t_{k+1} - t_k \equiv h, \; \forall k = -N_\delta, \ldots, -1.
\end{equation}
We then consider the discretized version of \eqref{def:Xt}--\eqref{def:cost}:
\begin{align}
    & X(t_{k+1}) = X(t_k) + b(t_k, X_{t_k}, \pi(t_k)) h + \sigma(t_k, X_{t_k}, \pi(t_k)) \Delta W(t_k), \label{def:Xtdiscrete} \\
    & \inf_{\{\pi(t_k)\}_{k=0}^{N_T-1}} \EE\left[\sum_{k=0}^{N_T-1} f(t_k, X_{t_k}, \pi(t_k)) h + g(X_T)\right], \label{def:costdiscrete}
\end{align}
where $X_{t_k}$ represents the path with $N_\delta$ lags, and $\Delta W(t_k)$ is the increment in Brownian motions:
\begin{equation}
    X_{t_k} = (X(t_{k-N_\delta}),\ldots, X(t_k)), \quad \Delta W(t_k) = W(t_{k+1}) - W(t_k).
\end{equation}

Regarding the discretized system \eqref{def:Xtdiscrete}, one would expect results similar to \cite[Section 4]{kushner2008numerical}. That is, as the mesh size $h \to 0$, the value function associated to \eqref{def:Xtdiscrete}--\eqref{def:costdiscrete} converges to the one of the original problem \eqref{def:Xt}--\eqref{def:cost}; and the near-optimal control associated to \eqref{def:Xtdiscrete} (which now as functions of $X_{t_k}$) is also near-optimal to the original one. Below we propose two architectures in deep learning for approximating $\pi(t_k) \in \RR^m$.

\subsection{Feedforward neural network}
A feedforward neural network is a composition of several fully-connected layers $F_{d_1, d_2}(x)$:
\begin{equation}
    F_{d_1, d_2}(x) = \rho(Ax + b) : \RR^{d_1} \to \RR^{d_2},
\end{equation}
where $x \in \RR^{d_1}$ and $F_{d_1, d_2}(x) \in \RR^{d_2}$ are the input and output of this layer, $A \in \RR^{d_2\times d_1} $ and $b \in \RR^{d_2}$ are the weight matrix and bias vector, and $\rho(\cdot)$ is the activation function applied on each element of the vector individually. Common choices of the activation function include rectified linear unit (ReLU), identity, sigmoid, hyperbolic tangent:
\begin{equation}
    \rho_{\text{ReLU}}(x) = \max\{x,0\}, \quad \rho_{\text{Id}}(x) = x, \quad \rho_{\text{s}}(x) = \frac{1}{1+e^{-x}}, \quad \rho_{\tanh}(x) = \tanh(x).
\end{equation}

Motivated by the path-dependent structure of the considered problems (the change of current state only depends on the history up to lag $\delta$), a natural idea is to approximate $\pi(t_k)$ by a feedforward neural network taking the state history up to lag $\bar{\delta}$ as the input. Note here it could be $\bar{\delta}\neq \delta$ since we may not know the underlying true $\delta$ a priori.
Without loss of generality, we assume $\bar{\delta}=N_{\bar{\delta}}h~ (N_{\bar{\delta}}\in\mathbb{N}^+$) and define
$\bar{X}_{t_k} \equiv (X(t_{k-N_{\bar{\delta}}}),\ldots, X(t_k), t_k)\in \RR^{n\times (N_{\bar{\delta}}+1)+1}$. Then we represent the policy as 
\begin{equation}\label{eq:ffXt}
    \pi(t_k) \approx F_{d_{I}, m} \circ F_{d_{I-1}, d_{I}} \cdots  F_{d_1, d_2}\circ F_{n\times (N_{\bar{\delta}}+1)+1, d_1}(\bar{X}_{t_k}),
\end{equation}
where $I$ is the number of hidden layers. In this case, the algorithm will produce feedback controls, \emph{i.e.}, controls that are adapted to the canonical filtration of $X$, denoted by $\mc{F}_t^X$. Also, with a fixed input dimension and the added time variable in the input, we are able to share the parameters of sub-neural networks, thus reducing the parameter number by a factor of $N_T$ compared to $N_T$ different networks at each timestamp.

\subsection{Recurrent neural network}\label{sec:lstm}
The idea of recurrent neural networks (RNNs)~\cite{rumelhart1986learning} is to make use of sequential information. 
They have shown great success in natural language processing, handwriting recognition, etc.~\cite{graves2013generating,graves2013speech,graves2009offline} The most common RNN is long short-term memory (LSTM)~\cite{hochreiter1997long}. The advantage of an LSTM is the ability to deal with the vanishing gradient problem and data with lags of unknown duration. 

An LSTM is composed of a series of units, each of which corresponds to a timestamp, and each unit consists of a cell and three gates: input gate, output gate, and forget gate. Among these components, the cell keeps track of the information received so far, the input gate captures to which extent new input information flows into the cell, the forget gate captures to which extent the existing information remains in the cell, and the output gate controls to which extent the information in the cell will be used to compute the output of the unit. 
In our case, the $k^{th}$ unit is responsible for approximating $\pi(t_k)$:
\begin{equation}
\label{eq:lstm_gate}
\begin{aligned}
    &\text{forget gate: } f_k = \rho_{\text{s}}(W_f x_k + U_f h_{k-1} + b_f),\\
    &\text{input gate: } i_k =  \rho_{\text{s}}(W_i x_k + U_i h_{k-1} + b_i), \\
    &\text{ontput gate: } o_k =  \rho_{\text{s}}(W_o x_k + U_o h_{k-1} + b_o), \\
    &\text{cell: } c_k = f_k \odot c_{k-1} + i_k \odot \rho_{\tanh}(W_c x_k + U_c h_{k-1} + b_c), \\
    &\text{output of the } k^{th} \text{ unit: } h_k = o_k \odot  \rho_{\tanh}(c_k),
\end{aligned}
\end{equation}
where the operator $\odot$ denotes the Hadamard product, $x_k$ denotes the $k^{th}$ input, $c_k$ stores the cell information, and $h_k, f_k, i_k, o_k$ and $c_k$ are all $d_h$-dimensional vectors.
We take $(X(t_0), t_0), (X(t_1), t_1), (X(t_2),t_2), \dots$ as the input sequence $x_0, x_1, x_2, \dots$ in practice, and specify the initial information of $h_{k}, c_{k}$ according to the discretized initial condition $\varphi$ (see details in Section~\ref{sec_implementation}). Then we take an affine transformation of $h_k$ as the proxy of $\pi(t_k)$:
\begin{equation}\label{eq:LSTM}
    \pi(t_k) \approx Wh_k + b.
\end{equation}

Although for both schemes \eqref{eq:LSTM} and \eqref{eq:ffXt}, the input dimensions keep constant as $k$ changes, using \eqref{eq:ffXt} requires prior knowledge of $\delta$. That is, for \eqref{eq:ffXt} which we feed the discretized state values \eqref{def:Xtdiscrete} of length $N_{\bar\delta} + 1$, to obtain the best performance, one needs to get an good estimate $\bar\delta$ of $\delta$ first;  while for \eqref{eq:LSTM} we only need to provide the current state value $X(t_k)$. Notice that in an LSTM all input information up to time $t_k$ is summarized by the $k^{th}$ cell, but if the optimal control were only depend on the past up to $\delta$, the forget gates are designed for dropping out the unneeded information. This dropout is characterized by NN's parameters, which are determined by supervised learning. We shall detail the learning part in the next section.

We remark that there are many variations of LSTM, for instance, gated recurrent units (GRUs)~\cite{cho2014learning} that do not have output gates, peephole LSTM~\cite{gers2002learning} where $h_{k-1}$ is mostly replaced by $c_{k-1}$ in all gates, etc. The numerical experiments will be conducted using the standard LSTM introduced above, and extensions to the variants are straightforward.

\subsection{Choice of input data: $X(t_k)$ or $W(t_k)$}
In the previously proposed networks, we use the data consisting of the state $X(t_k)$ as the input.
On one side, using $X(t_k)$ as the input data may lead to sub-optimal controls, as $(\mc{F}_t^X)_{t\geq 0}$ might be smaller than $(\mc{F}_t)_{t\geq 0}$ in general, and thus the resulting policy in the feedback form forms a strict subset of \eqref{def_A}. On the other side, in many scenarios, there exists an optimal control in \eqref{def_A} of the feedback form. For instance, see examples of problems with delay, among many others, in \cite{chen2010maximum,chen2012delayed,oksendal2000maximum}. This is also the case in control problems without delay, if the solution to the HJB equation is smooth enough, which provides the decoupling field of the corresponding forward-backward stochastic differential equations \cite[Chapter 4]{carmona2016lectures}. So we do not lose much by searching within a smaller set.

A naive idea to enforce the progressive measurability imposed in \eqref{def_A} is to take the data consisting of the Brownian motion $\{W(t_k)\}_{k=1}^{N_T}$ as input.
However, in numerical experiments, we observed that the resulting networks become much more challenging to optimize, leading to much worse performance. One possible explanation is that when there are additional constraints on the admissible set $\mathbb{A}$, it is usually directly related to the state process $X(t)$. The constraints are hardly satisfied if NNs only know background noises $W$ and need to infer the state $X$.
Another drawback of using  $W(t_k)$ as the input data in the feedforward model is that the whole input variable becomes $W_k=(W(t_0), \cdots, W(t_k))$, whose dimension increases dramatically as $k$ approaches to $N_T$. In addition, due to different input dimensions across $k$, we cannot let sub-neural networks share parameters, which will further increase the number of parameters.
Based on the above reasons, we only report results with the input data consisting of $X(t_k)$.

\subsection{Implementation}
\label{sec_implementation}
To summarize the algorithms proposed in the above subsections, we essentially have
\begin{equation}
    \pi(t_k) = \psi_k(\text{Data} ; \theta),
\end{equation}
where Data can be either $(X(t_{k-N_{\bar{\delta}}}),\ldots, X(t_k), t_k)$ in the feedforward model or $(X(t_k), t_k)$ in the LSTM model, and $\theta$ denotes all parameters appearing in \eqref{eq:ffXt} or \eqref{eq:LSTM}.

The internal states $h_k, c_k$ in the LSTM model~\eqref{eq:lstm_gate} need to be initialized properly according to the initial segment $\varphi$. In our case, we start the LSTM model from the timestamp $t_{-N_\delta}=-\delta$ (when $\delta=\infty$, we choose a properly truncated time), and the initial states $h_{-N_{\delta}}$ and $c_{-N_{\delta}}$ are both initialized as $(X(-\delta = t_{-N_\delta}), 0,\cdots,0)$. Here we have assumed $d_h\geq n$ (recall that $n$ is the dimension of $X(t)$) such that there is no information lost at the beginning, and additional zeros are added to match the dimension $d_h$. Then we feed in the input sequence $(X(t_{-N_{\delta}+1}), t_{-N_{\delta}+1}), (X(t_{-N_{\delta}+2}), t_{-N_{\delta}+2}), \dots$ to evolve the model according to \eqref{eq:lstm_gate} and output controls starting at $t_0$.

The optimal parameters $\theta^\ast$ are then be obtained by minimizing the following expected discretized loss using stochastic gradient descent algorithms:
\begin{equation}
    \inf_{\{\psi_k \in \mc{N}_k\}_{k=0}^{N_T-1}}  \EE\left[\sum_{k=0}^{N_T-1} f(t_k, X_{t_k}, \psi_k(\text{Data}; \theta)) h + g(X_T)\right].
\end{equation}
Note that in some numerical examples, we may need to deal with constraints involving the control policy and/or state variables. When the policy $\pi$ taking values in $\mc{A}$ is required to be nonnegative, we apply a ReLU activation function before the final output of the policy network to ensure such a property:
\begin{equation*}
    \pi(t_k) \leftarrow \rho_{\text{ReLU}}(\pi(t_k)).
\end{equation*}
When a 1-dimensional state variable $X(t)$ is required to be nonnegative, we add the corresponding penalty term in the cost function:
\begin{align*}
    & f(t_k, X_{t_k}, \psi_k(\text{Data}; \theta)) \leftarrow f(t_k, X_{t_k}, \psi_k(\text{Data}; \theta)) + \lambda \max\{-X(t_{k+1}),0\}
\end{align*}
where $\lambda$ is a penalty coefficient. More complicated constraints can be dealt with by the penalty method in a similar way (see \emph{e.g.},~\cite{HaE:16}).
Details on the choices of $\psi_k$, algorithm to obtain $\theta^\ast$, $N_T$, etc., are presented at the beginning of Section~\ref{sec_numerics}.


\section{Benchmark examples}\label{sec_benchmark}
This section presents three tractable examples: a linear-quadratic regulator problem with delay in engineering, an optimal consumption problem in a financial market with delayed dynamics, and a portfolio optimization problem with infinite delay. They together serve as preparation of numerical experiments in Section~\ref{sec_numerics}.
With specific formulas of $(b, \sigma, f, g)$ depending on the current state $X(t)$, the weighted average of $X_t$ and $X(t-\delta)$, the first two problems turn out to be essentially finite-dimensional and admit solutions expressed in a simple form. This allows us to benchmark and compare our proposed two deep learning schemes. A third example, where the portfolio performance depends on the exponential average of all the historical value ($\delta = \infty$), is presented also with analytical solutions, to further evident the superiority of the LSTM model. In the sequel, we will focus on describing the model and keep technical details minimal. The optimal control $\pi^\ast$ and cost $V_0$ to the problem \eqref{def:Xt}--\eqref{def:cost}:
\begin{equation}
V_0 := \sup_{\pi \in \mathbb{A}}\EE_{\varphi}\left[\int_t^T f(s, X_s, \pi(s)) \ud s  + g(X_T) \Big\vert X_0 = \varphi \right]
\end{equation}
are provided in Propositions~\ref{prop:LQ}--\ref{prop:portfolio} and we give references on their proofs.

For the first two examples, we take a special dependence form for the functionals $(b, \sigma, f, g)$. To be specific, the path dependence is characterized by distributed delay $Y(t)$ and discrete delay $Z(t)$:
\begin{equation}\label{def_dependent}
  Y(t) := \int_{-\delta}^0 e^{\lambda s} X(t+s) \ud s \text{ and } Z(t) := X(t-\delta). 
\end{equation}
Note that $Y(t)$ is also called the exponentially decayed weighted moving average. 
Problems with such explicit structures have been studied in \cite{oksendal2000maximum,elsanosi2000some,elsanosi2001optimal,larssen2001hjb,bauer2005stochastic,oksendal2011optimal}. Some allow general analysis and some have tractable examples. We choose two examples from \cite{bauer2005stochastic} and present them in Sections~\ref{sec:LQ} and \ref{sec:consumption}.

\subsection{Linear-quadratic problem with delay}\label{sec:LQ}
Stochastic linear-quadratic (LQ) problems were extensively studied in the literature. They have appeared in many contexts and have been used to benchmark various numerical algorithms due to their tractability. LQ problems with delay was first investigated by Kolmanovski{\u{\i}} and Sha{\u{\i}}khet \cite{kolmanovskiui1996control}. The delay version can be stated as: 
\begin{equation}\label{eq:lq}
(LQ)~ 
\begin{dcases}
\ud X(t) = (A_1(t)X(t) + A_2(t)Y(t) + A_3 Z(t) + B(t)\pi(t)) \ud t + \sigma(t) \ud W(t), \quad t  \in[0,T]\\
\min_{\pi(t) \in \RR^m}\EE_{\varphi} \left[\int_0^T (X(t) + e^{\lambda \delta}A_3 Y(t))\transpose Q(t) (X(t) + e^{\lambda \delta}A_3 Y(t)) + \pi(t)\transpose R(t)\pi(t)\ud t \right.\\
\hspace{75pt} + (X(T) + e^{\lambda\delta} A_3Y(T))\transpose G (X(T) + e^{\lambda\delta} A_3Y(T))\Bigg],
\end{dcases}
\end{equation}
where $X_0 = \varphi \in L^2(\Omega, C)$ is a given initial segment, $A_1(t), A_2(t), Q(t) \in \RR^{n\times n}$, $B(t) \in \RR^{n \times m}$, $R(t) \in \RR^{m\times m}$ are deterministic matrix-valued functions in $L^\infty[0,T]$, $\sigma(t) \in \RR^{n \times \ell}$ is a deterministic matrix-valued function in $L^2[0,T]$, $A_3, G \in \RR^{n \times n}$ are deterministic matrices. We assume that $Q(t), G$ are positive semi-definite and $R(t)$ is positive definite for all $t \in [0,T]$ and continuous on $[0,T]$. To have a tractable solution, we further prescribe the relation:
\begin{equation}\label{lq:parameters}
A_2(t) = e^{\lambda\delta}(\lambda I_n + A_1(t) + e^{\lambda\delta} A_3) A_3,
\end{equation}
where $I_n$ is the identity matrix with rank $n$. This example was studied in
\cite[Section 4]{bauer2005stochastic}, and we summarize the main results as follows for completeness.
\begin{prop}\label{prop:LQ}
	Consider the stochastic control problem (LQ) with the initial segment $\varphi$. Assume that $A_3 \neq 0$ and \eqref{lq:parameters} holds, then the optimal control is given by
	\begin{equation}
	\pi^\ast(t) = -R\inv(t) B(t)\transpose P(t)(X^\ast(t) + e^{\lambda\delta}A_3Y^\ast(t)), \text{ where } Y^\ast(t) = \int_{-\delta}^0 e^{\lambda s} X^\ast(t+s) \ud s, 
	\end{equation}
	$X^\ast(t)$ solves the SDDE \eqref{eq:lq} with the optimal control $\pi^\ast$, and $P(t)$ solves the Riccati equation
	\begin{equation}
	\dot P(t) = P(t) B(t) R\inv(t)B(t)\transpose P(t) - (A_1(t) + e^{\lambda\delta}A_3)\transpose P(t) - P(t)(A_1(t) + e^{\lambda\delta}A_3) - Q(t), \quad P(T) = G.
	\end{equation}
The optimal cost $V_0$ to problem (LQ)   is given by
	\begin{equation}
	V_0 \equiv V(0, X_0) = (X(0) + e^{\lambda\delta}A_3Y(0))\transpose P(0) (X(0) + e^{\lambda\delta}A_3Y(0)) + \int_0^T \emph{Tr}(\sigma(s)\sigma\transpose(s) P(s))\ud s.
	\end{equation}
\end{prop}

\subsection{Optimal consumption in a delayed financial market}\label{sec:consumption}
In this section, we consider $X(t)$ described in the SDDE \eqref{def:Xt} as the wealth process, and study the utility maximization problem from both consumption and terminal wealth. Here we do not specify dynamics of tradable assets, but directly model the wealth return and volatility, both depending on the sliding average $Y(t) = \int_{-\delta}^0 e^{\lambda s} X(t+s) \ud s$ and the past value $Z(t) = X(t-\delta)$. This can be interpreted as investing on of some path-dependent options. Let $c(t)$ be the investor's consumption rate at time $t$, then the optimal consumption problem is stated as: 
\begin{equation}\label{eq:consumption}
(C)~
\begin{dcases}
\ud X(t) = (\mu(t, X(t), Y(t)) + aZ(t) - c(t)) \ud t + \sigma(t, X(t), Y(t)) \ud W(t), \quad t  \in[0,T]\\
\max_{c(t) \in \RR^+} \EE_{\varphi} \left[\int_0^T e^{-\beta t} U_1(c(t)) \ud t + U_2(X(T) + ae^{\lambda \delta}Y(T)) \right], \text{ subject to } X(t) \geq 0, \text{ } t \in [0, T],\\
\end{dcases}
\end{equation}
where $U_1, U_2$ are utility functions on the consumption and the terminal wealth, and $a$ is a positive real number. In this example, all processes are scalars, \emph{i.e.}, $n = m = \ell = 1$. To have a tractable solution, we take the power utility $U_1(x) = U_2(x) = \frac{1}{\gamma}x^\gamma$ and require 
\begin{equation}\label{consumption:parameters}
\mu(t,x,y) = ae^{\lambda\delta}(ae^{\lambda\delta} + \lambda)y + \mu_t T(x,y), \quad
\sigma(t,x,y) = \sigma_t T(x,y),
\end{equation}
where $T(x,y) = x + ae^{\lambda\delta}y$, and $\mu_t$ and $\sigma_t$ are some positive  continuous functions. This example has been studied in \cite[Section 5]{bauer2005stochastic}, and a special case $a = -\lambda e^{-\lambda \delta}$, $\lambda <0$ was treated in \cite{elsanosi2001optimal}. We now summarize the results as follows.

\begin{prop}\label{prop:consumption}
	Consider the optimal consumption problem (C) with the initial wealth segment $X_0 = \varphi$, and assume \eqref{consumption:parameters} holds. Then, the problem value $V_0$ is given by
	\begin{equation}
	V_0 \equiv V(0, X_0) =\frac{1}{\gamma} p(0)^{1-\gamma} (X(0) + ae^{\lambda\delta}Y(0))^\gamma,
	\end{equation}
	where  $p(t)$ solves the following equation:
	\begin{equation}
	\dot p(t) =(\half \gamma \sigma_t^2 - \frac{\gamma}{1-\gamma}(\mu_t + ae^{\lambda\delta}))\cdot p(t) - e^{-\frac{\beta t}{1-\gamma}}, \quad p(T) = 1. 
	\end{equation}
    The optimal consumption $c^\ast(t)$ is
	\begin{equation}
	c^\ast(t) = e^{-\frac{\beta t}{1-\gamma}} \frac{1}{p(t)} (X^\ast(t) + ae^{\lambda\delta}Y^\ast(t)),
	\end{equation}
	with $X^\ast(t)$ following the SDDE \eqref{eq:consumption} associated with $c^\ast$.
\end{prop}

\begin{rem}
If the terminal utility is changed to $e^{-\beta_2T} U_2(X(T) + ae^{\lambda\delta}Y(T))$, then Proposition~\ref{prop:consumption} still holds except for the terminal condition of $p(t)$. It then becomes
\begin{equation}
p(T) = e^{-\beta_2T/(1-\gamma)}.
\end{equation}
If the running utility is changed to $U_1(x) = \frac{\eta x^\gamma}{\gamma}$, the ordinary differential equation for $p(t)$ then becomes
\begin{equation}
\dot p(t) =(\half \gamma \sigma_t^2 - \frac{\gamma}{1-\gamma}(\mu_t + ae^{\lambda\delta}))\cdot p(t) - \eta^{\frac{1}{1-\gamma}}e^{-\frac{\beta t}{1-\gamma}}.
\end{equation}
\end{rem}

\subsection{Portfolio optimization with complete memory}
\label{sec:portfolio}
In the last example, we study a portfolio optimization problem with an infinite delay feature. More precisely, in the SDDE \eqref{def:Xt} and cost functional \eqref{def:cost}, the dependence is on $X(t)$ and the exponential average of all the history value $Y(t)$,
\begin{equation}\label{def:Yt}
Y(t) := \int_{-\infty}^0 e^{\lambda s}X(t+s)\ud s.
\end{equation}
Let $X(t)$ be the wealth process as in Section~\ref{sec:consumption}. Here we consider both investment $\pi(t)$ and consumption $c(t)$ and parameterize them proportional to the wealth $X(t)$, \emph{i.e.}, $c(t)$ denotes the fraction of wealth consumed at time $t$. Then the problem reads:
\begin{equation}\label{eq:portfolio}
(P)~
\begin{dcases}
\ud X(t) = [((\mu_1 - r)\pi(t) - c(t) + r)X(t) + \mu_2 Y(t)] \ud t + \sigma\pi(t)X(t) \ud W(t), \quad t  \in[0,T]\\
\max_{\substack{c(t) \in \RR^+ \\ \pi(t) \in \RR}} \EE_{\varphi} \left[\int_0^T e^{-\beta t} U_1(c(t)X(t)) \ud t + e^{-\beta T} U_2(X(T),Y(T)) \right].
\end{dcases}
\end{equation}
The above problem was introduced in \cite{pang:17}. They derived explicit solutions under exponential, power and log utilities. Given the history up to time 0, $\varphi(\cdot) \in L^2(\Omega, C((-\infty, 0], \RR))$, the problem value turns out depending only on $X(0) \equiv \varphi(0)$ and $Y(0)  \equiv \int_{-\infty}^0 e^{\lambda s} \varphi(t + s) \ud s$.

\begin{rem}
For different choices of utility functions, we may or may not require $X(t) \geq 0$, for $t \in [0,T]$. Theoretically, if the strategy $(\pi(t), c(t))$ satisfies 
\begin{align}
     \abs{\pi(t)X(t)} \leq \Lambda_0 \abs{X(t) + Y(t)}, \\
    \abs{c(t)X(t)} \leq \Lambda_0 \abs{X(t) + Y(t)},
\end{align}
for some constant $\Lambda_0$, then $X(t)$ in \eqref{eq:portfolio} stays positive. For proof, see \cite[Lemma 2.1]{pang:17}.
\end{rem}


\begin{prop}\label{prop:portfolio}
    Consider the portfolio optimization problem  with infinite delay (P) with the complete history $X_0 = \varphi \in L^2(\Omega, C((-\infty, 0], \RR))$ under log utility: 
    \begin{align}
        U_1(x) = \log(x),  \quad U_2(x, y) = \frac{1}{\beta} \log (x+ \eta y),  \quad \eta = \half\left(\sqrt{(r+\lambda)^2 + 4\mu_2} - (r+\lambda) \right).
    \end{align}
    and require a state constraint $X(t) \geq 0$, $\forall t \in [0, T]$. Then the optimal investment and consumption rates are
    \begin{align}
        &\pi^\ast(t) = \frac{(\mu_1 - r)(X^\ast(t) + \eta Y^\ast(t))}{\sigma^2 X^\ast(t)}, \\
        &c^\ast(t) = \frac{\beta (X^\ast(t) + \eta Y^\ast(t))}{X^\ast(t)},
    \end{align}
    where $X^\ast(t)$ solves the SDDE \eqref{eq:portfolio} associated to $\pi^\ast$ and $c^\ast$. The value function is of the form
    \begin{equation}
        V(t,x,y) = p(t) + \frac{1}{\beta}\log(x+ \eta y)
    \end{equation}
    with 
    \begin{equation}
        p(t) = \frac{\Lambda_2}{\beta}(1- e^{-\beta(T-t)}), \quad \Lambda_2 = \half \frac{(\mu_1-r)^2}{\beta \sigma^2} + \log(\beta) - 1 + \frac{1}{\beta}(r + \eta),
    \end{equation}
    and the optimal cost $V_0 = V(0, X(0), Y(0))$.
\end{prop}

\section{Numerical results}\label{sec_numerics}
In this section, we present numerical results on the three benchmark examples introduced above: linear-quadratic problem~(Section~\ref{sec_num_lq}), optimal consumption~(Section~\ref{sec_num_csmp}), and portfolio optimization~(Section~\ref{sec_num_polog}). We choose $n = 10$ in the first example, while $n = 1$ is fixed in the settings of the other two examples. In all three experiments, we observe consistently good performance compared to the benchmark solutions. The important hyperparameters and running time of these problems are summarized in Table~\ref{tab_parameter}, and other problem-dependent parameters will be introduced in the corresponding subsections. 
The initial path $\varphi \in L^2(\Omega, C)$ is modeled as a fixed deterministic path perturbed by some white noise.

In numerical solutions, the distributed delay $Y(t)$ is approximated by the midpoint quadrature rule, in both calculating the benchmark solution and simulating the dynamics of $X(t)$ (cf. equations \eqref{eq:lq},\eqref{eq:consumption} and \eqref{eq:portfolio}). The activation function in the feedforward model is ReLU in all the layers except that in the final output layer, it is an identity if there is no constraint on the control.
In the learning process, to guarantee a fair comparison, in each problem, we use the same learning rate and choose the size of the feedforward network and the LSTM properly such that they have a similar number of parameters.
We adopt Adam optimizer~\cite{Kingma2015adam} to optimize the parameters in neural networks. In each update step, we simulate 128 (batch size) paths to compute the stochastic gradient through backpropagation. 
For each example, we run the algorithm three times and report the average result.

The algorithm is implemented in Python by using the machine learning library TensorFlow~\cite{abadi2016tensorflow}. The code can be found in a public GitHub repository\footnote{\url{https://github.com/frankhan91/RNN-ControlwithDelay}} upon publication, and thus the results presented here can be straightforwardly reproduced and further developed.
\begin{table}[h]
	\caption{Some hyperparameters and running time of the numerical examples}\label{tab_parameter}
	\begin{center}
		\begin{tabular}{@{}c|ccc@{}}
			\toprule
			Parameter / Problem & Linear-Quadratic  & Optimal Consumption & Portfolio Optimization\\ \midrule
			$T$         &  1 &  2 &  5  \\
			$\delta$    &  1 &  0.8 &  $~~\infty~^{\dagger}$  \\
			$N_T$       &  40 &  60 &  100  \\
			$\xi$ (constraint penalty)   &  0 (no constraint) &  $10^5$ &  10 \\
			running time (s) & 4294 & 1429 & 941 \\
			\bottomrule
		\end{tabular}
	\end{center}
	\small{$^\dagger$The initial segment $\varphi$ is constant when $t\leq -4$ to ease the computation of the distributed delay $Y(t)$. Both neural network models are initialized at $t=-4$.
	}
\end{table}

\subsection{Linear-quadratic problem}
\label{sec_num_lq}
\begin{figure}[!htb]
\centering
\includegraphics[width=0.5\textwidth]{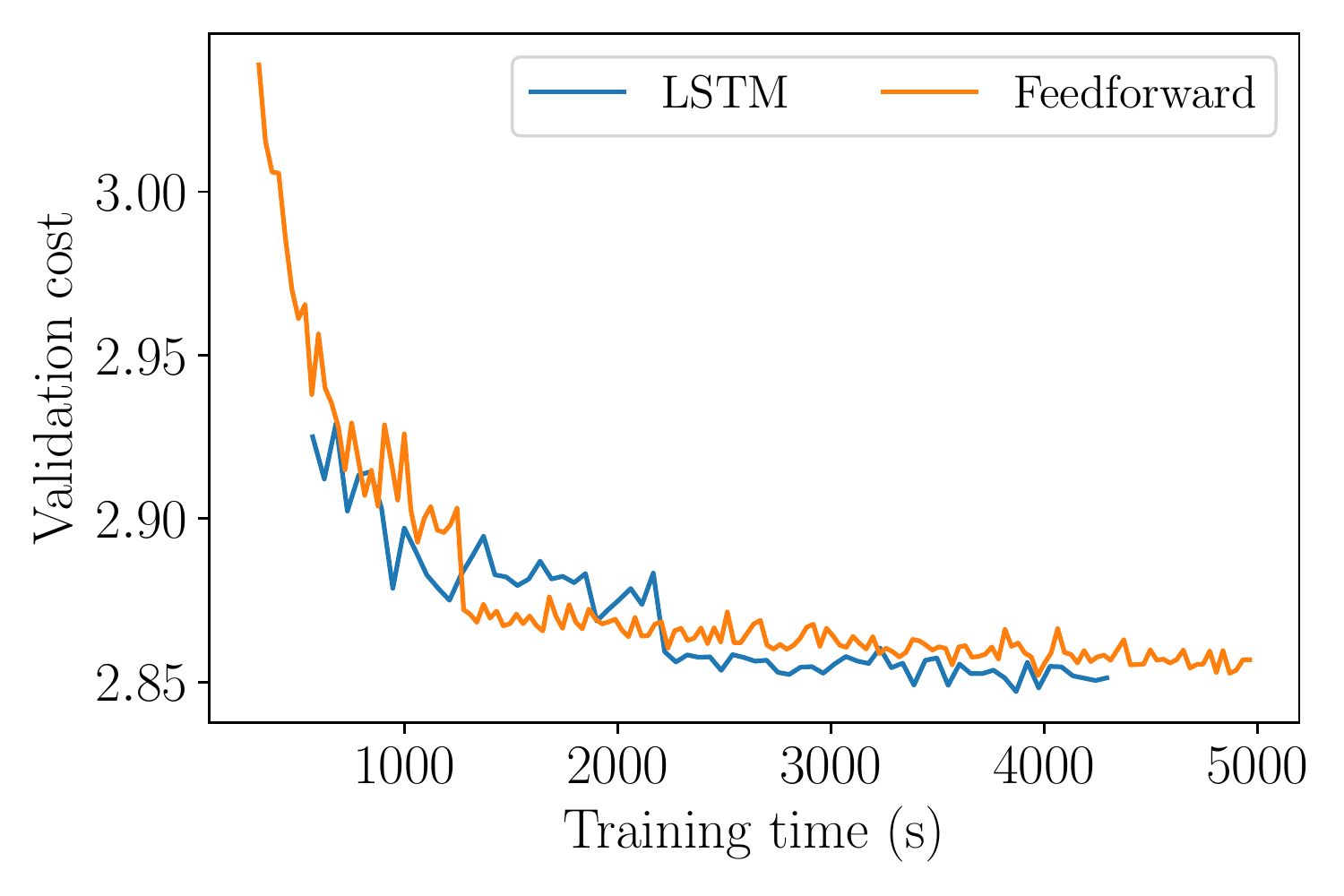}
  \caption{Training curve of two models in the example of linear-quadratic problem.
  }
  \label{fig:lq_train_curve}
\end{figure}

\begin{figure}[!htb]
\centering
\includegraphics[width=0.98\textwidth]{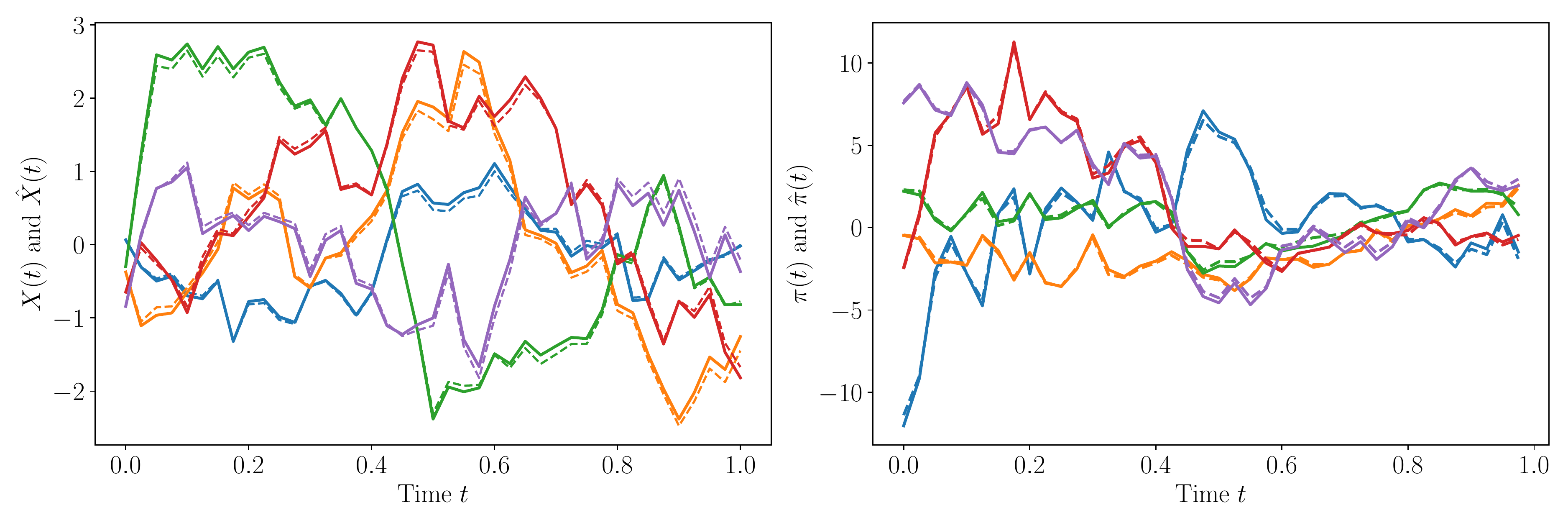}
  \caption{
      A sample path of the first 5 dimensions of the state $X(t)$ and control $\pi(t)$ obtained from the LSTM model.
      Left: the optimal state process discretized from the analytical solution $X_i(t)$ (solid lines) and its approximation $\hat{X}_i(t)$ (dashed lines) provided by the approximating control, under the same realized path of Brownian motion.
      Right: comparisons of the optimal control $\pi^*_i(t)$ (solid lines) and $\hat{\pi}_i(t)$ (dashed lines).
  }
  \label{fig:lq_path_lstm}
\end{figure}

\begin{figure}[!htb]
\centering
\includegraphics[width=0.98\textwidth]{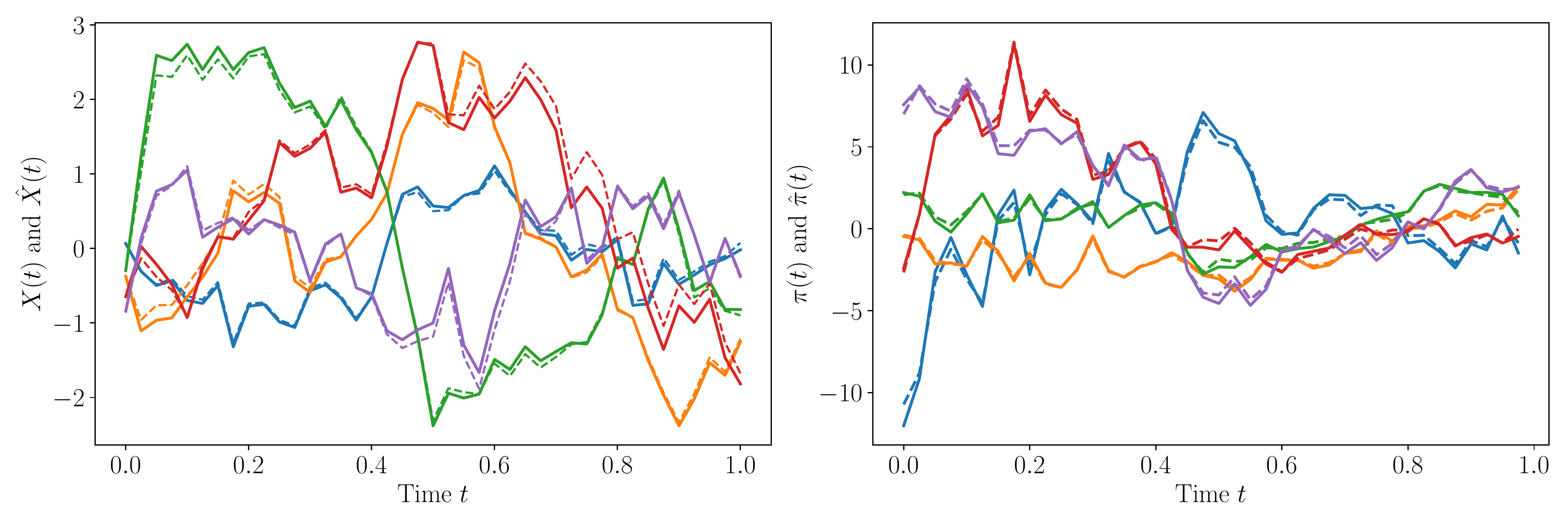}
  \caption{
      A sample path of the first 5 dimensions of the state $X(t)$ and control $\pi(t)$ obtained from the feedforward model.
      Left: the optimal state process $X_i(t)$ discretized from the analytical solution (solid lines) and its approximation $\hat{X}_i(t)$ (dashed lines) provided by the approximating control, under the same realized path of Brownian motion.
      Right: comparisons of the optimal control $\pi^*_i(t)$ (solid lines) and $\hat{\pi}_i(t)$ (dashed lines).
  }
  \label{fig:lq_path_shff}
\end{figure}

We consider a 10-dimensional example in Section~\ref{sec:LQ}, in which $n=10$, $m=10$, and $\lambda=0.1$.
In Eq.~\eqref{eq:lq}, $A_1, A_3, B, \sigma$ are random generated constant coefficient matrices, $Q, R, G$ are constant matrices proportional to identity matrices, and $A_2$ is determined by~\eqref{lq:parameters}.
The dimension of the hidden state in the LSTM model is $d_h = 200$, which gives 171610 parameters in total. The feedforward model takes the state history as inputs up to lag $\bar{\delta}=\delta$ with $N_{\bar{\delta}}=40$, and it has 2 hidden layers with a width of 300, which gives 108910 parameters in total. The LSTM model is trained with 16000 steps, in which the learning rate is 0.005 for the first 8000 steps and 0.0005 for the second 8000 steps. The feedforward model is trained with 32000 steps, in which the learning rate is 0.005 for the first 8000 steps and 0.0005 for the remaining 24000 steps. We calculate the total cost on the validation data every 200 steps and plot the curve against training time in Fig.~\ref{fig:lq_train_curve}. We can see that with the same learning rates and roughly the same number of parameters, the LSTM model converges to a solution with a lower cost compared to the feedforward model. We further test the learned policy on a new set of sample paths, and the average cost of the LSTM model and feedforward model are 2.8740 and 2.8812. Note that the cost obtained by the policy discretized from the analytical optimal control is 2.8723 (cf. Proposition~\ref{prop:LQ}).
Figs.~\ref{fig:lq_path_lstm} and \ref{fig:lq_path_shff} depict one sample path (first 5 dimensions only) of the optimal state $X$ and control $\pi$ provided by two neural networks in comparison with the analytical solution, in which the LSTM architecture presents a better agreement.

\begin{figure}[!htb]
\centering
\includegraphics[width=0.5\textwidth]{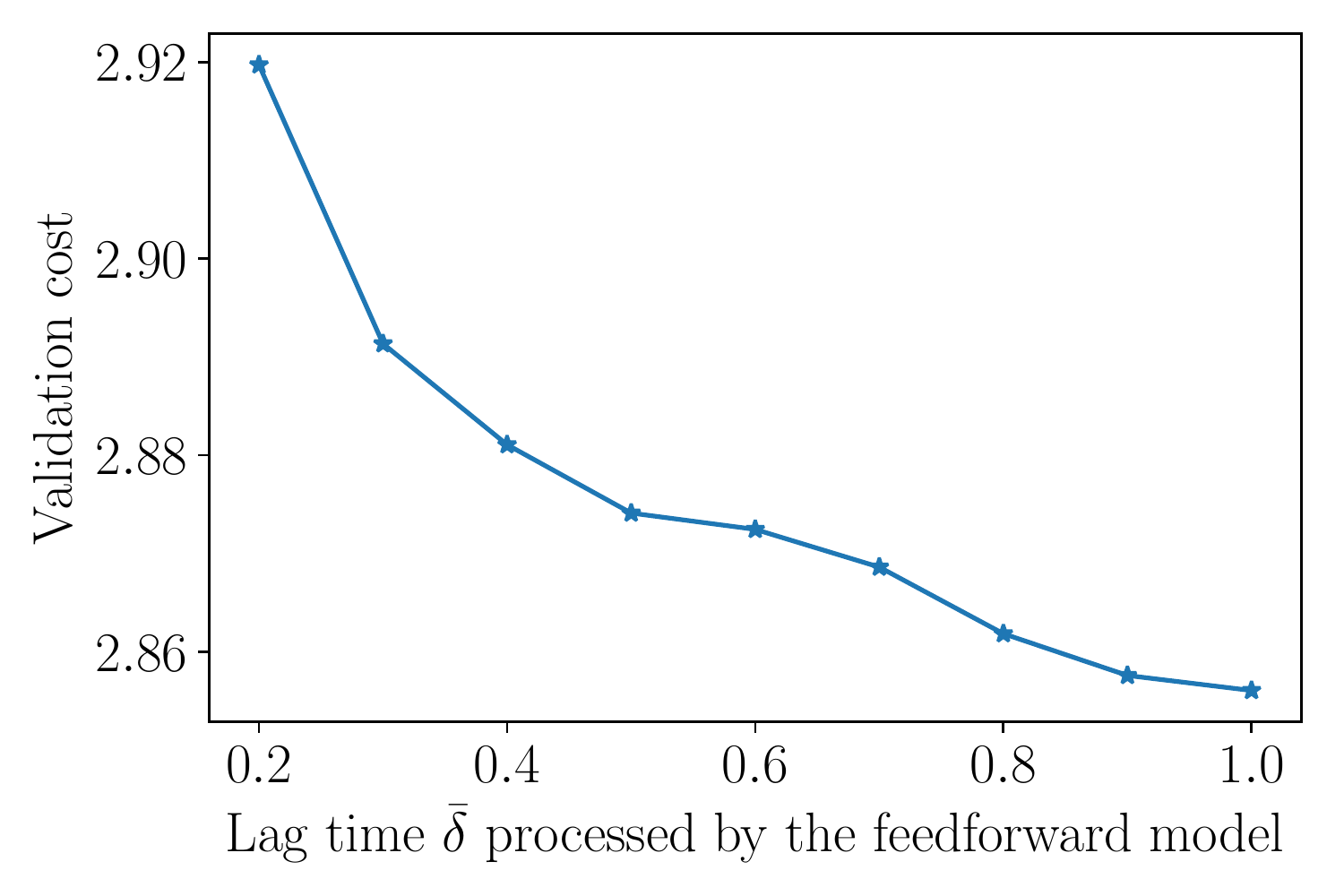}
  \caption{The effect of lag time $\bar{\delta}$ processed by the feedforward model in the example of linear-quadratic problem. The lag time $\delta$ in the actual system is 1.}
  \label{fig:lq_shff_lag}
\end{figure}

In the above experiment, the lag time $\bar{\delta}$ processed by the feedforward model is chosen to be the same as $\delta$.
One main drawback of the feedforward model is that it requires to know the true lag time $\delta$ a priori to determine the network's size. If the chosen lag time $\bar{\delta}$ is smaller than the actual lag $\delta$ time, which means there is a loss of information when the feedforward network processes the data, the final performance might be compromised. To quantify this effect, we test the feedforward model with different processed lag time $\bar{\delta}$ from $0.2$ to $1$ with step size 0.1, while the actual lag $\delta=1$. The corresponding optimized costs are shown in Fig.~\ref{fig:lq_shff_lag}. As expected, we observe that the cost increases as the lag time processed by the feedforward model decreases. A higher optimized cost indicates that the model can only find a sub-optimal but not an optimal strategy due to the lack of information. 

\subsection{Optimal consumption}
\label{sec_num_csmp}
\begin{figure}[!htb]
\centering
\includegraphics[width=0.5\textwidth]{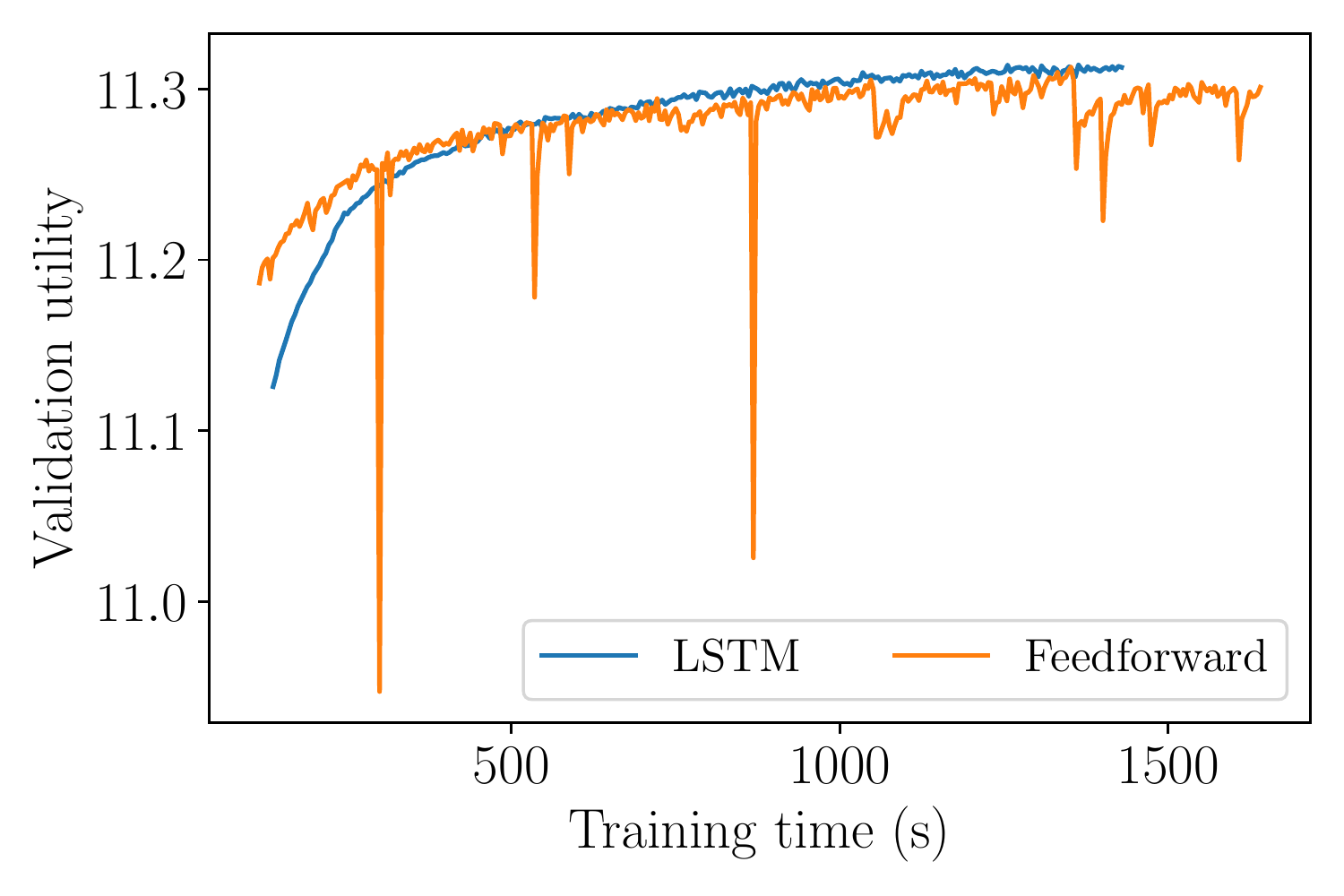}
  \caption{Training curve of two models in the example of optimal consumption.
  }
  \label{fig:csmp_train_curve}
\end{figure}

\begin{figure}[!htb]
\centering
\includegraphics[width=0.98\textwidth]{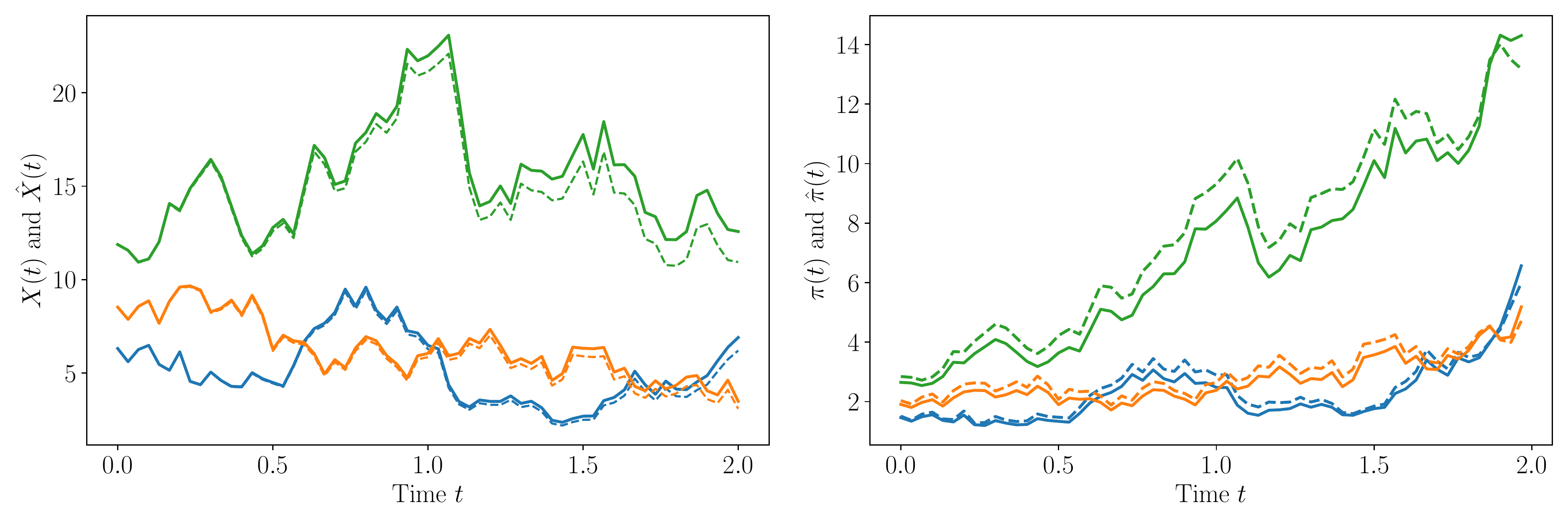}
  \caption{
      Three sample paths of the state $X(t)$ and control $\pi(t)$ obtained from the LSTM model.
      Left: the optimal state process $X(t)$ discretized from the analytical solution (solid lines) and its approximation $\hat{X}(t)$ (dashed lines) provided by the approximating control, under the same realized paths of Brownian motion.
      Right: comparisons of the optimal control $\pi^*_i(t)$ (solid lines) and $\hat{\pi}_i(t)$ (dashed lines).
  }
  \label{fig:csmp_path_lstm}
\end{figure}

\begin{figure}[!htb]
\centering
\includegraphics[width=0.98\textwidth]{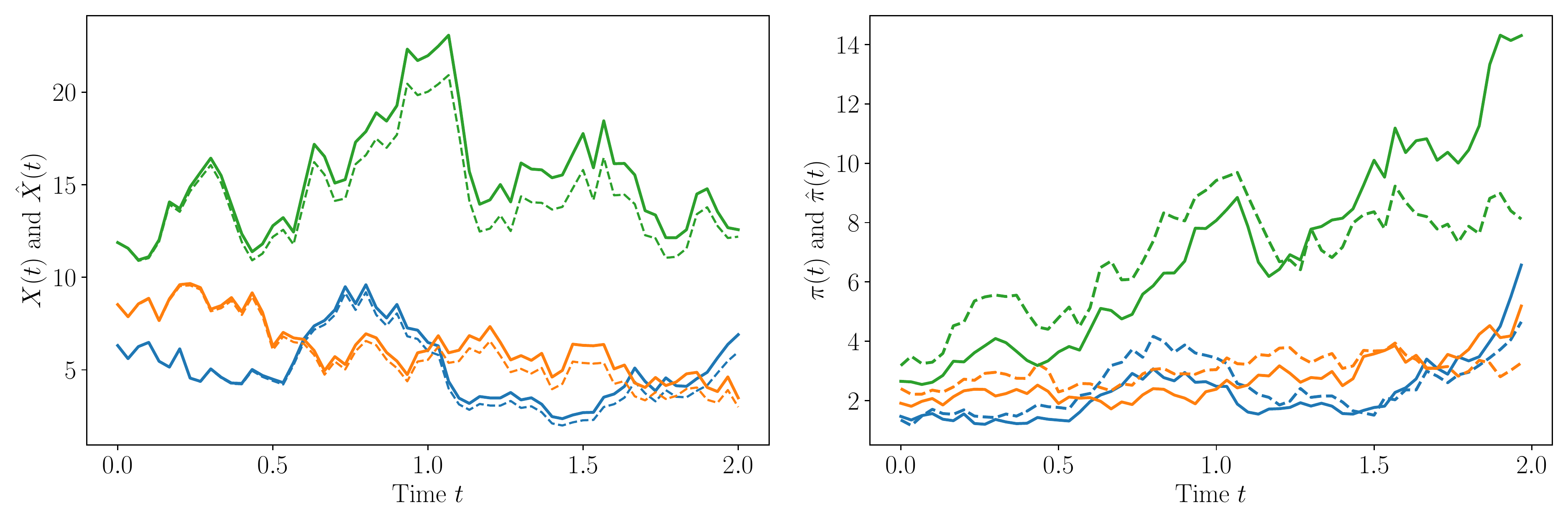}
  \caption{
      Three sample paths of the state $X(t)$ and control $\pi(t)$ obtained from the feedforward model.
      Left: the optimal state process $X(t)$ discretized from the analytical solution (solid lines) and its approximation $\hat{X}(t)$ (dashed lines) provided by the approximating control, under the same realized paths of Brownian motion.
      Right: comparisons of the optimal control $\pi^*_i(t)$ (solid lines) and $\hat{\pi}_i(t)$ (dashed lines).
  }
  \label{fig:csmp_path_shff}
\end{figure}

We consider the example introduced in Section~\ref{sec:consumption}, with $\lambda=0.1$, $\beta=0.1$, $a=0.2$, $\gamma=0.7$, $\mu_t=0.1$, and $\sigma_t=0.5$.
The dimension of the hidden state in the LSTM model is $d_h = 30$, which gives 3991 parameters in total. The feedforward model takes the state history as input up to lag $\bar{\delta}=\delta$ with $N_{\bar{\delta}}=24$, and it has 2 hidden layers with a width of 54, which gives 4483 parameters in total. 
The LSTM model and feedforward model are trained with 60000 and 80000 steps, respectively, with a learning rate of 0.0001.
We calculate the total utility on the validation data every 200 steps and plot the curve against training time in Fig.~\ref{fig:csmp_train_curve}. We can see that with the same learning rates and a similar number of parameters, the LSTM model still converges to a higher utility (which corresponds to a better control policy) than the feedforward model. Furthermore, the training process of the feedforward model is much more unstable than the LSTM model, probably due to its insufficient capability to handle state constraints. We further test the learned policy on a new set of sample paths, and the average utility of the LSTM model and feedforward model are 10.7202 and 10.7100. Note that the utility obtained by the policy discretized from the analytical optimal control is 10.7142 (cf. Proposition~\ref{prop:consumption}).
Figs.~\ref{fig:csmp_path_lstm} and \ref{fig:csmp_path_shff} depict three sample paths provided by two models in comparison with the analytical solution, in which the LSTM model presents a better agreement.

\subsection{Portfolio optimization}
\label{sec_num_polog}
\begin{figure}[!htb]
\centering
\includegraphics[width=0.5\textwidth]{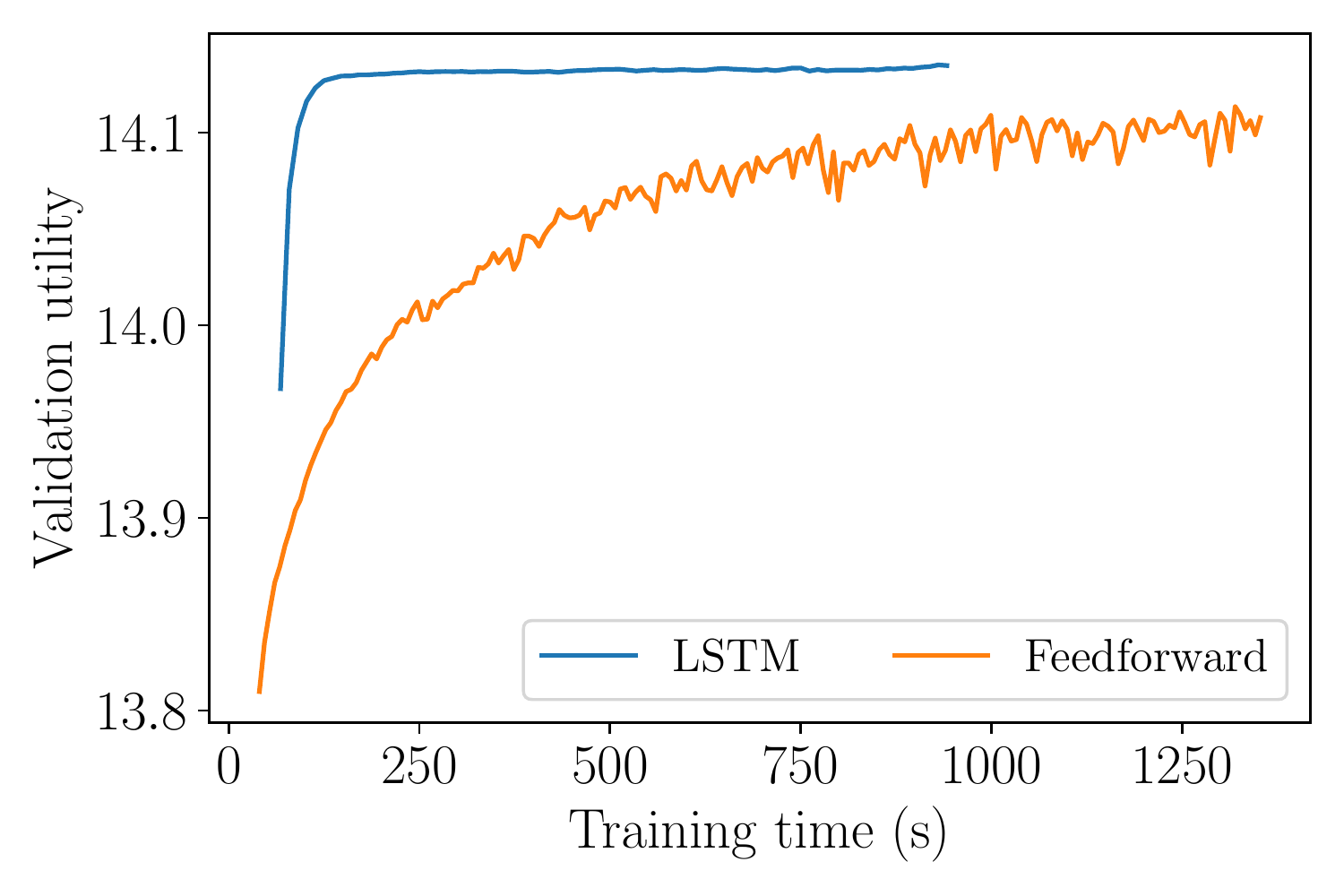}
  \caption{Training curve of two models in the example of portfolio optimization.
  }
  \label{fig:polog_train_curve}
\end{figure}

\begin{figure}[!htb]
\centering
\includegraphics[width=0.98\textwidth]{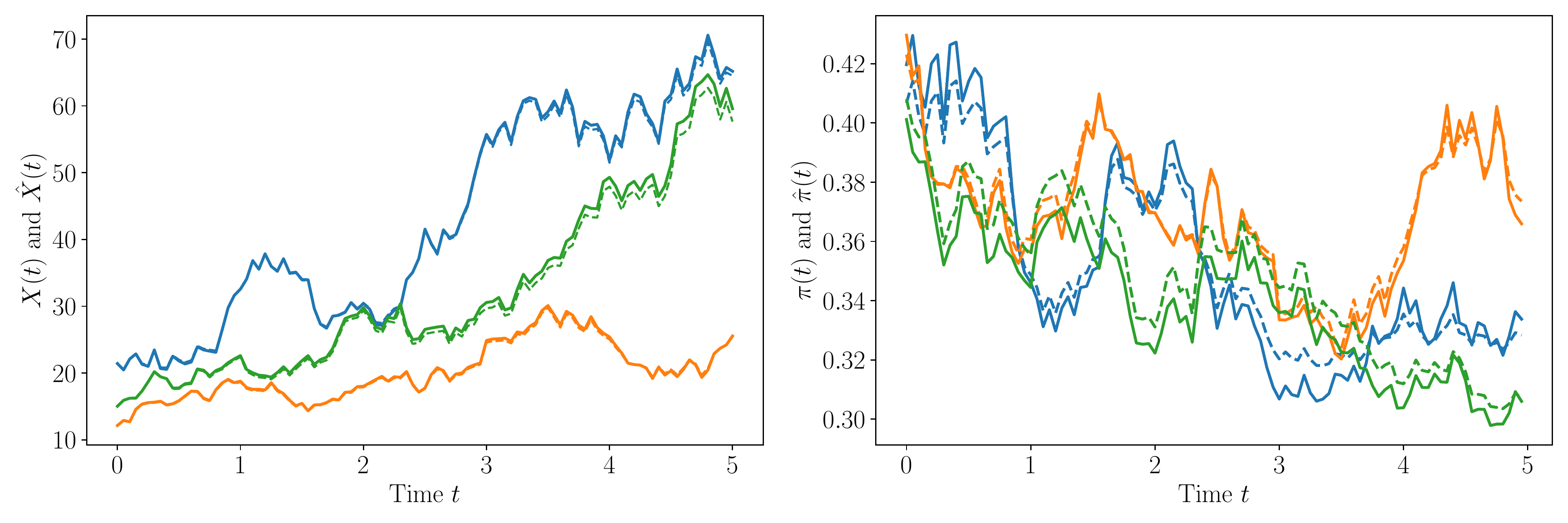}
  \caption{
      Three sample paths of the state $X(t)$ and control $\pi(t)$ obtained from the LSTM model.
      Left: the optimal state process $X(t)$ discretized from the analytical solution (solid lines) and its approximation $\hat{X}(t)$ (dashed lines) provided by the approximating control, under the same realized paths of Brownian motion.
      Right: comparisons of the optimal control $\pi^*_i(t)$ (solid lines) and $\hat{\pi}_i(t)$ (dashed lines).
  }
  \label{fig:polog_path_lstm}
\end{figure}

\begin{figure}[!htb]
\centering
\includegraphics[width=0.98\textwidth]{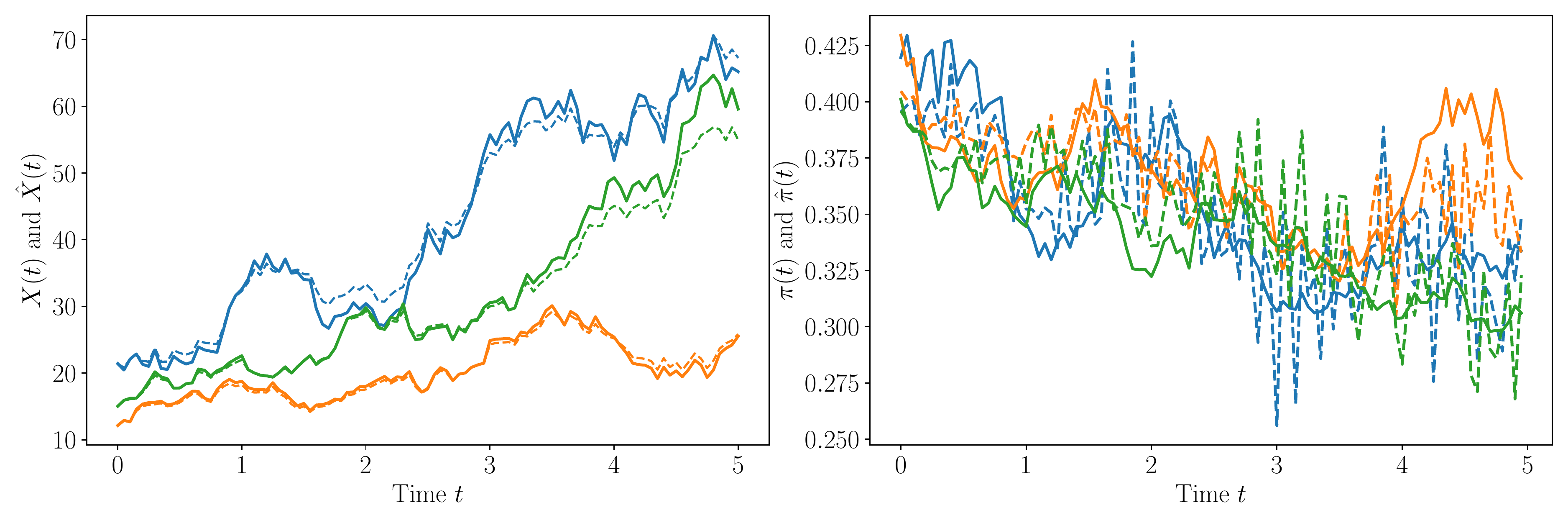}
  \caption{
      Three sample paths of the state $X(t)$ and control $\pi(t)$ obtained from the feedforward model.
      Left: the optimal state process $X(t)$ discretized from the analytical solution (solid lines) and its approximation $\hat{X}(t)$ (dashed lines) provided by the approximating control, under the same realized paths of Brownian motion.
      Right: comparisons of the optimal control $\pi^*_i(t)$ (solid lines) and $\hat{\pi}_i(t)$ (dashed lines).
  }
  \label{fig:polog_path_shff}
\end{figure}

We now consider the example introduced in Section~\ref{sec:portfolio}, with $\lambda=0.3$, $\beta=0.2$, $\mu_1=0.1$, $\mu_2=0.2$, $r=0.05$, and $\sigma=0.4$. In this example, the model has a infinite memory. We choose the initial segment $\varphi$ to be constant for $t\leq -4$ and the white noise perturbation only presents for $t\in[-4,0]$. Under this setting, the computation of the distributed delay $Y(t)$, defined as an integral running from $-\infty$ to $t$, becomes tractable.
Both neural network models are initialized at $t=-4$ for the sake of fairness.
The dimension of the hidden state in the LSTM model is $d_h = 60$, which gives 15242 parameters in total. The feedforward model takes the state history as input up to lag $\bar{\delta}=4 < \delta=\infty$ with $N_{\bar{\delta}}=80$, and it has 2 hidden layers with a width of 96, which gives 17474 parameters in total.
The LSTM model and the feedforward model are trained with 8000 and 20000 steps, respectively, both with a learning rate of $5\times10^{-5}$.
We calculate the total utility in \eqref{eq:portfolio} on the validation data every 100 steps and plot the curve against the training time in Fig.~\ref{fig:polog_train_curve}. We can see that in this example, the LSTM model converges to a higher utility much faster than the feedforward model, and the convergence is much stable. We further test the learned policy on a new set of sample paths, and the average utilities obtained from the LSTM model and feedforward model are 14.5326 and 14.4979. Note that the utility obtained by the policy discretized from the analytic optimal control is 14.5316 (cf. Proposition~\ref{prop:portfolio}). Fig.~\ref{fig:polog_path_lstm} and \ref{fig:polog_path_shff} depict three sample paths provided by two models in comparison with the analytical solution, in which the LSTM model presents a much better agreement. 

We believe that the remarkable outperformance of the LSTM model is due to the infinite memory feature of this problem, which can be naturally captured by the LSTM model but not the feedforward one. 
In particular, for the feedforward model, the first input at $t=0$ is $(X(-4), \ldots, X(0), 0)$, and later inputs keep the same sliding window length ($N_{\bar{\delta}}=80$) of the state history and gradually drop the initial information. We believe that such information loss is the main reason why the feedforward model performs significantly sub-optimal. Numerically, we also observed that the choice of $\bar\delta < 4$ will even worsen the performance due to the information loss at the initial step.

\section{Conclusion}\label{sec_conclusion}
In this paper, we propose and systematically study deep neural networks-based algorithms to solve stochastic control problems with delay features. The challenge is brought by the path-dependent property in the SDDE system and thus its intrinsic high dimensions. Viewing the optimal policy as a function of state process path $(X(t))_{t \geq 0}$ or background noise path $(W(t))_{t \geq 0}$, we naturally start with feedforward neural networks and take the discretized process as inputs. With the path-dependent feature, we then employ recurrent neural networks for sequence modeling to parameterize the policy and optimize the objective function. On the architecture design level, we point out that recurrent neural networks such as the LSTM model do not require prior knowledge of lag time $\delta$. The models are then tested on three benchmark examples: 1. linear-quadratic problem with fixed delay; 2. optimal consumption in a financial market with fixed finite delay; 3. portfolio optimization with complete memory. Numerically, we found that the architecture of recurrent neural networks can naturally capture the path-dependent feature with much flexibility, resulting in a better performance with more efficient and stable training compared to the feedforward architectures. The superiority is even evident for infinite delay $\delta = \infty$, supported by our third example. Moreover, the carefully selected benchmark problems with open-sourced code will facilitate the further study of numerical algorithms for stochastic control problems with delay features. Remark that the numerical techniques developed in this paper can be naturally generalized to previous studies \cite{Hu2:19, HaHu:19, HaHuLo:20} on the deep fictitious play algorithm for finding the Nash equilibrium in stochastic differential games by including delay in the games. This will be investigated in future work.

\bibliographystyle{plain}
\bibliography{Reference}

\end{document}